\documentclass[11pt,reqno]{amsart}
\usepackage{amssymb,latexsym,mathrsfs,
	amsxtra,amscd,amsfonts,amsthm,amsmath,verbatim,epsfig}

\usepackage{color}
\usepackage[colorlinks,pagebackref=true]{hyperref}

\usepackage[english]{babel}

\makeatletter

\@addtoreset{equation}{section}
\def\theequation{\thesection.\@arabic \c@equation}
\makeatletter
\def\widebreve#1{\mathop{\vbox{\m@th\ialign{##\crcr\noalign{\kern3\p@}%
				\brevefill\crcr\noalign{\kern3\p@\nointerlineskip}%
				$\hfil\displaystyle{#1}\hfil$\crcr}}}\limits}
\def\brevefill{$\m@th \setbox\z@\hbox{$\braceld$}%
	\bracelu\leaders\vrule \@height\ht\z@ \@depth\z@\hfill\braceru$}

\makeatletter

\def\@citecolor{blue}
\def\@linkcolor{red}
\def\@urlcolor{blue}
\def\@urlcolor{blue}\def\theequation{\arabic{equation}}
\def\theequation{\thesection.\arabic{equation}}
\numberwithin{equation}{section}

\def\dim{\operatorname{dim}}
\def\depth{\operatorname{depth}}
\def\height{\operatorname{ht}}
\def\grade{\operatorname{grade}}
\def\ass{\operatorname{Ass}}
\def\supp{\operatorname{Supp}}

\def\reg{\operatorname{reg}}

\def\spec{\operatorname{Spec}}

\def\ZZ{\mathbb Z}
\def\QQ{\mathbb Q}

\newcommand{\NN}{\mathbb N}
\newcommand{\mm}{\mathfrak m}

\newcommand{\ov}{\overline}
\newcommand{\lm}{{\lambda}}

\newcommand{\R}{\mathcal R}

\newcommand{\n}{{\bf{n}}}

\newcommand{\m}{{\bf {m}}}

\newcommand{\J}{J(1)_{1}\cdots J(r)_{1}}

\newcommand{\bl}{\begin{lemma}}
	\newcommand{\el}{\end{lemma}}
\newcommand{\bt}{\begin{theorem}}
	\newcommand{\et}{\end{theorem}}
\newcommand{\ben}{\begin{enumerate}}
	\newcommand{\een}{\end{enumerate}}
\newcommand{\bpf}{\begin{proof}}
	\newcommand{\eepf}{\end{proof}}
\newcommand{\beqn}{\begin{eqnarray*}}
	\newcommand{\eeqn}{\end{eqnarray*}}
\newcommand{\beqnn}{\begin{eqnarray}}
\newcommand{\eeqnn}{\end{eqnarray}}
\newcommand{\bd}{\begin{definition}}
	\newcommand{\ed}{\end{definition}}
\newcommand{\bp}{\begin{proposition}}
	\newcommand{\ep}{\end{proposition}}
\newcommand{\bc}{\begin{corollary}}
	\newcommand{\ec}{\end{corollary}}
\newcommand{\bex}{\begin{example}}
	\newcommand{\eex}{\end{example}}

\newcommand{\wlg}{ Without loss of generality }

\theoremstyle{plain}

\newtheorem{theorem}{Theorem}[section]
\newtheorem{corollary}[theorem]{Corollary}

\newtheorem{proposition}[theorem]{Proposition}

\newtheorem{lemma}[theorem]{Lemma}

\newtheorem{example}[theorem]{Example}
\newtheorem{definition}[theorem]{Definition}

\newtheorem{question}[theorem]{Question}

\theoremstyle{remark}
\newtheorem{remark}[theorem]{Remark}
\numberwithin{equation}{theorem}

\begin{document}
	\title[Rees' theorem for filtrations, multiplicity function and reduction criteria]{Rees' theorem for filtrations, multiplicity function and reduction criteria}
	\author{Parangama Sarkar}
	\address{chennai mathematical institute, h1, sipcot it park, siruseri, chennai, india-603103}
	\email{parangamasarkar@gmail.com,  parangama@cmi.ac.in}
	
	\thanks{The author is supported by Department of Science \& Technology (DST); INSPIRE Faculty Fellowship.} 
	
	\begin{abstract}
	Let $J\subset I$ be ideals in a formally equidimensional local ring with $\lm(I/J)<\infty.$ Rees proved that for all $n\gg 0$, $\lm(I^n/J^n)$ is a polynomial $P(I/J)(X)$ in $n$ of degree at most $\dim R$ and  $J$ is a reduction of $I$ if and only if $\deg P(I/J)(X)\leq \dim R-1.$ We extend this result for all Noetherian filtrations of ideals in a formally equidimensional local ring  and for (not necessarily Noetherian) filtrations of ideals in analytically irreducible rings.  We provide certain classes of ideals such that $\deg P(I/J)$ achieves its maximal degree. On the other hand, for ideals $J\subset I$ in a formally equidimensional local ring, we consider the multiplicity function $e(I^n/J^n)$ which is a polynomial in $n$ for all large $n.$ We explicitly determine the $\deg e(I^n/J^n)$ in some special cases. For  an ideal $J$ of analytic deviation one, we give characterization of reductions in terms of $\deg e(I^n/J^n)$ under some additional conditions.    
	\end{abstract}
	\maketitle
	
	\thispagestyle{empty}
\section{introduction}
	Let $(R,\mm)$ be a Noetherian local ring of dimension $d$ and $I$ be an ideal in $R.$ If $I$ is $\mm$-primary then Samuel \cite{P1} showed that for all large $n,$ the {\it Hilbert-Samuel function of $I,$} $\displaystyle H_I(n)=\lambda({R}/{I^n})$ (here $\lambda(M)$ denotes the length of an $R$-module $M$) coincides with a polynomial $$\displaystyle P_I(n)=e_0(I)\binom{n+d-1}{d}-e_1(I)
\binom{n+d-2}{d-1}+\cdots+(-1)^d e_d(I)$$ of degree $d,$ called the {\it Hilbert-Samuel polynomial of $I.$} The coefficient $e(I):=e_0(I)$ is called the {\it multiplicity of $I.$} Rees used this numerical invariant $e(\_)$ to study the numerical characterization of reductions. He showed that if $R$ is a formally equidimensional local ring (i.e. $R$ is a Noetherian local ring and its completion in the topology defined by the maximal ideal is equidimensional) and $J\subset I$ are $\mm$-primary ideals in $R$ then $J$ is a reduction of $I$ if and only if $e(I)=e(J)$ \cite{Rees4}. In literature this result has been generalized in several directions. Amao \cite{Am}, considered more general case where $I,J$ are not necessarily $\mm$-primary ideals and proved that if $M$ is a finitely generated $R$-module and $J\subset I$ are ideals of $R$ such that $IM/JM$ has finite length then $\mu(n):=\lm(I^nM/J^nM)$ is a polynomial $P(IM/JM)$ in $n$ for $n\gg 0.$ In \cite{Rees}, Rees showed that in a Noetherian local ring $(R,\mm),$ the degree of $P(I/J)$ is at most $\dim R$ and he further proved the following theorem.
\bt{\em(Rees, \cite[Theorem 2.1]{Rees})}
	Let $R$ be a formally equidimensional Noetherian local ring and $J\subset I$ be ideals in $R$ with $\lambda(I/J)<\infty.$ Then $J$ is a reduction of $I$ if and only if  $\deg P(I/J)$ is at most $\dim R-1.$
\et
Since for some $n,$ $J^n$ is a reduction of $I^n$ implies $J$ is a reduction of $I,$ we can interpret Rees' result as follows:
\bt\label{Reestheorem}
	 	Let $R$ be a formally equidimensional Noetherian local ring of dimension $d$ and $J\subset I$ be ideals in $R$ with $\lambda(I/J)<\infty.$  Then the following are true.
\ben
	\item  The function $F(n)=\lim_{m\to\infty}\lambda(I^{mn}/J^{mn})/m^d$ is a polynomial in $n$ of degree at most $d.$ 
	\item Fix $n\in\mathbb Z_+,$ then $J$ is a reduction of $I$ if and only if $F(n)=0.$\een
\et

In \cite[Theorem 6.1]{C2014}, Cutkosky proved the following result.
\bt{\em(Cutkosky, \cite[Theorem 6.1]{C2014})}
	Let $R$ be an analytically unramified local ring of dimension $d>0.$ Suppose $\{I_i\}$ and $\{J_i\}$ are two graded families of ideals such that for all $i,$ $J_i\subset I_i$ and there exists $c\in\mathbb Z_+,$ such that $\mm^{ci}\cap I_i=\mm^{ci}\cap J_i$ for all $i.$ Assume that if $P$ is a minimal prime of $R$ then $I_1\subset P$ implies $I_i\subset P$ for all $i>0.$ Then the limit $\lim_{m\to\infty}\lambda(I_{m}/J_{m})/m^d$ exists.  
\et
Therefore it is natural to ask whether Rees' result (Theorem \ref{Reestheorem}) holds true  for two (not necessarily Noetherian) filtrations of ideals $\{I_i\}$ and $\{J_i\}$ in an analytically unramified local ring $R$ of dimension $d.$ In this direction we first prove Rees' theorem for Noetherian filtrations of ideals. 
\bt{\em(=Theorem \ref{Noetheriancase})}
	Let $(R,\mm)$ be a  Noetherian local ring of dimension $d>0$ and for $1\leq l\leq r,$ $\mathcal I(l)=\{I(l)_n\},$ $\mathcal J(l)=\{J(l)_n\}$ be Noetherian filtrations  of ideals in $R$ such that $J(l)_n\subset I(l)_n$ and $\lambda( I(l)_n/J(l)_n)<\infty$ for all $n\in\mathbb N$ and $1\leq l\leq r.$ Let $$R_1=\bigoplus_{m_1,\ldots,m_r\geq 0}I(1)_{m_1}\cdots I(r)_{m_r}\mbox{ and } R_2=\bigoplus_{m_1,\ldots,m_r\geq 0}J(1)_{m_1}\cdots J(r)_{m_r}.$$ Fix $n_1,\ldots,n_r\in\mathbb Z_+.$ If $R_1$ is integral over $R_2$ then
	$$\lim_{m\to\infty}\lambda(I(1)_{mn_1}\cdots I(r)_{mn_r}/J(1)_{mn_1}\cdots J(r)_{mn_r})/m^d=0.$$
Converse holds if $R$ is formally equidimensional and $\grade(\J)\geq 1.$
	\et
We prove that the ``only if" part of Rees' theorem \ref{Reestheorem}$(2)$  holds true for (not necessarily Noetherian) filtrations of ideals in analytically irreducible local rings (i.e. $R$ is a Noetherian local ring and its completion in the topology defined by the maximal ideal is domain). We also give an example to show that the ``if" part of Rees' theorem \ref{Reestheorem}$(2)$ is not true for non-Noetherian filtrations of ideals in general (See example \ref{counter}).
\bt{\em(=Theorem \ref{domain})}
	Let $(R,\mm)$ be  an analytically irreducible local ring of dimension $d>0$ and for $1\leq l\leq r,$ $\mathcal I(l)=\{I(l)_n\},$ $\mathcal J(l)=\{J(l)_n\}$ be (not necessarily Noetherian) filtrations of ideals in $R$ such that $J(l)_n\subset I(l)_n$ for all $n\in\mathbb N$ and $1\leq l\leq r.$ Fix $n_1,\ldots,n_r\in\mathbb Z_+.$ Suppose there exists an integer $c\in\mathbb Z_+$ such that for all $i\in\mathbb N,$ $$\mm^{ci}\cap I(1)_{in_1}\cdots I(r)_{in_r}=\mm^{ci}\cap J(1)_{in_1}\cdots J(r)_{in_r}. $$ Let $R_1=\bigoplus_{m_1,\ldots,m_r\geq 0}I(1)_{m_1}\cdots I(r)_{m_r}$ be integral over $R_2=\bigoplus_{m_1,\ldots,m_r\geq 0}J(1)_{m_1}\cdots J(r)_{m_r}.$ Then $$\lim_{m\to\infty}\lambda(I(1)_{mn_1}\cdots I(r)_{mn_r}/J(1)_{mn_1}\cdots J(r)_{mn_r})/m^d=0.$$
\et
Rees' theorem \cite[Theorem 2.1]{Rees} shows that $J$ is a reduction of $I$ implies  $P(I/J)$ is a polynomial of degree at most $\dim R-1.$  In this direction one can ask that if $J$ is a reduction of $I$ when the degree of the polynomial $P(I/J)$ attains the maximal. We use the notation $(I,J)$ to denote $J$ is a reduction of $I.$ We  produce a lower bound of degree of $ P(I/J)$ and show that if $J$ is a complete intersection ideal then $\deg P(I/J)=l(J)-1$ where $l(J)$ is the analytic spread of $J$ (See Theorem \ref{n4}). We also improve the lower bound of $\deg P(I/J)$ for the following class of ideals. 
\bt{\em(=Theorem \ref{gra})}
	Let $(R,\mm)$ be a Noetherian local ring, $J\subsetneq I$ be ideals in $R$ such that $\lm(I/J)<\infty$ and $\grade (J)=s.$ Let $a_1,\ldots,a_s,b_1,\ldots,b_t$ be a $d$-sequence in $R$ which minimally generates the ideal $J$ and $t\geq 1.$  Suppose $I\cap(J':b_1)=J'$ where $J'=(a_1,\ldots,a_s)$ $(J'=(0)$ if $s=0)$. Then the following are true.
	\ben
	\item $\grade (J)\leq \deg P(I/J).$
	\item If $J$ is a reduction of $I$ then $\grade (J)\leq \deg P(I/J)\leq \l(J)-1.$ 
	\een
\et
This helps us to provide a large class of ideals $(I,J)$ which attain the maximal degree of $ P(I/J),$ i.e. $l(J)-1.$ 
\bc{\em(=Corollary \ref{Artin})}
		Let $(R,\mm)$ be a Cohen-Macaulay local ring of dimension $d\geq 2$ with infinite residue field, $I$ be an ideal of $R$  and $J$ be any  minimal reduction of $I$ with $\lm(I/J)<\infty.$  Suppose $I$ satisfies $G_{l(I)}$ and $AN_{l(I)-2}^{-}$ conditions. Then for any ideal $K$ with  $J\subsetneq K\subseteq I,$  $\deg P(K/J)=l(J)-1.$
	\ec
In the other direction, Herzog, Puthenpurakal and  Verma  \cite{HPV} proved that, for ideals $J\subset I$ in a formally equidimensional local ring, $e(I^n/J^n)$ is eventually a polynomial function. In the same vein, in a recent paper, Ciuperc\u{a} \cite{cc15}, showed that $e(I^n/J^n)$ is a polynomial function of degree at most $\dim R-t$ where $t$ is a constant equal to $\dim Q$ where $Q=R/J^n:I^n$ for all $n\gg 0$ and in \cite{cc16}, he remarked that if $J$ is reduction of $I$ then $\deg e(I^n/J^n)\leq l(J)-1.$ He also studied the function $e(I^n/J^n)$ where $J$ is an equimultiple ideal (i.e. $\height J=l(J)$) and proved the following characterization of reductions for equimultiple ideals in terms of $\deg e(I^n/J^n)$ \cite[Theorem 2.6]{cc16}.
\bt{\em(Ciuperc\u{a}, \cite[Theorem 2.6]{cc16})}\label{cata}
Let $(R,\mm)$ be a formally equidimensional local ring and $J\subseteq I$ proper ideals of $R$ with $J$ equimultiple. Let $f(n)=e(I^n/J^n).$ The following are true.
\ben
\item If $J\subseteq I$ is a reduction then $\deg f(n)\leq l(J)-1.$
\item If $J\subseteq I$ is not a reduction then $\deg f(n)= l(J).$
\een
\et
In this situation one may ask {\it if $J$ has analytic deviation one (i.e. $l(J)=\height(J)+1$) can we give characterization of reduction  in terms of $\deg e(I^n/J^n)?$}
\\Motivated by the  result of Ciuperc\u{a} (Theorem \ref{cata}), we provide upper and lower bounds of $\deg e(I^n/J^n)$ and show that if $J$ is a complete intersection ideal then $J$ is reduction of $I$ if and only if $\deg e(I^n/J^n)= l(J)-1$ (See Proposition \ref{d1}). In the case $J$ has analytic deviation one, we characterize when $J$ is a reduction of $I$ in terms of $\deg e(I^n/J^n)$ under some additional hypotheses.
\bt {\em(=Theorem \ref{ad1})}
	Let $(R,\mm)$ be a formally equidimensional local ring  of dimension $d\geq 2,$ $J\subsetneq I$ be ideals in $R$ and $J$ has analytic deviation one. Suppose $l(J_p)< l(J)$ for all prime ideals $p$ in $R$ such that $\height p=l(J).$ Then the following are true.
	\ben
	\item If $J$ is not a reduction of $I$ then $\deg e(I^n/J^n)=l(J)-1.$
	\item If $l(J)=d-1,$ $\depth (R/J)>0$ and for all $n\geq 1,$ $\sqrt{J:I}=\sqrt{J^n:I^n}$ then \\$J$ is a reduction of $I$ if and only if $\deg e(I^n/J^n)\leq l(J)-2.$
	\een
\et
We can not omit the condition $\depth (R/J)>0$ from Theorem \ref{ad1} $(2)$ (see Example \ref{ex2}). We give sufficient conditions on the ideal $J$ for the equality $\sqrt{J:I}=\sqrt{J^n:I^n}$ for all $n\geq 1$ (see  Proposition \ref{req}).
\section{Rees' theorem for filtrations}
In this section we prove Rees' theorem for (not necessarily Noetherian) filtrations of ideals. A family $\mathcal I=\{I_n\}_{n\in\mathbb N}$ of ideals in $R$ is called {\it filtration of ideals} if $I_0=R,$ $I_m\subseteq I_n$ for all $m\geq n$ and $I_mI_n\subseteq I_{m+n}$ for all $m,n\in\mathbb N.$ A filtration $\mathcal I=\{I_n\}_{n\in\mathbb N}$ of ideals in $R$ is said to be {\it Noetherian filtration} if  $\oplus_{n\geq 0} I_n$ is a finitely generated $R$-algebra. We first prove Rees' result for Noetherian filtrations of ideals and then using iterated Noetherian filtrations, we prove the ``only if" part for (not necessarily Noetherian) filtrations of ideals in analytically irreducible rings. We conclude this section with an example which shows that for non-Noetherian filtrations the ``if" part of Rees' theorem is not true  in general.
\bl\label{use}
	Let $(R,\mm)$ be a local ring of dimension $d>0$ and $A=\bigoplus_{m_1,\ldots,m_r\geq 0}A_{m_1,\ldots,m_r}$ be a $\mathbb N^r$-graded finitely generated $R$-algebra where $A_{0,\ldots,0}=R.$ Let $u_1,\ldots,u_r\in\mathbb Z_+.$ Consider the grdaed subring $A^{(u_1,\ldots,u_r)}=\bigoplus_{m_1,\ldots,m_r\geq 0}A_{u_1m_1,\ldots,u_rm_r}$ of $A.$ Then $A$ is finitely generated $A^{(u_1,\ldots,u_r)}$-module.
\el
\begin{proof}
	Let $\{a_1,\ldots,a_p\}$ be a generating set of $A$ as $R$ algebra. Define $\mathcal S=\{a_1^{r_1}\cdots a_p^{r_p}\mid 0\leq r_i< u_1\cdots u_r\}.$ Let $h_1,\ldots,h_p\in\mathbb N$ such that $h_i=u_1\cdots u_rq_i+t_i$ where $0\leq t_i<u_1\cdots u_r$ and $1\leq i\leq p.$ Then $a_1^{h_1}\cdots a_p^{h_p}=(a_1^{q_1}\cdots a_p^{q_p})^{u_1\cdots u_r}a_1^{t_1}\cdots a_p^{t_p}.$ Thus $\mathcal S$ is a generating set of $A$ as  $A^{(u_1,\ldots,u_r)}$-module.
\end{proof}
\bl\label{reduction}
Let $(R,\mm)$ be a Noetherian local ring and for $i=1,\ldots,s,$ $J_i\subset I_i$ be ideals in $R$ and $\grade (J_1\cdots J_s)\geq 1.$ Suppose $J_i$ is not a reduction of $I_i$ for some $i\in \{1,\ldots,s\}.$ Then $J_1\cdots J_s$ is not a reduction of $I_1\cdots I_s.$
\el
\begin{proof}
	Without loss of generality, assume that $J_1$ is not a reduction of $I_1.$ Suppose  for some $n\geq 0,$ $(J_1\cdots J_s)(I_1\cdots I_s)^n=(I_1\cdots I_s)^{n+1}.$ Then  $$J_1I^{n}_1(I_2\cdots I_s)^{n+1}\supset (J_1I^n_1)\cdots(J_sI^n_s)=I_1I^{n}_1(I_2\cdots I_s)^{n+1}\supset J_1I^{n}_1(I_2\cdots I_s)^{n+1}.$$ Hence $J_1M=I_1M$ where $M=I^{n}_1(I_2\cdots I_s)^{n+1}.$ Therefore by \cite[Lemma 1.5]{Rees2}, $J_1$ is a reduction of $I_1$ which is a contradiction.
\end{proof}
\bt\label{Noetheriancase}
Let $(R,\mm)$ be a  Noetherian local ring of dimension $d>0$ and for $1\leq l\leq r,$ $\mathcal I(l)=\{I(l)_n\},$ $\mathcal J(l)=\{J(l)_n\}$ be Noetherian filtrations  of ideals in $R$ such that $J(l)_n\subset I(l)_n$ and $\lambda( I(l)_n/J(l)_n)<\infty$ for all $n\in\mathbb N$ and $1\leq l\leq r.$ Let $$R_1=\bigoplus_{m_1,\ldots,m_r\geq 0}I(1)_{m_1}\cdots I(r)_{m_r}\mbox{ and } R_2=\bigoplus_{m_1,\ldots,m_r\geq 0}J(1)_{m_1}\cdots J(r)_{m_r}.$$ Fix $n_1,\ldots,n_r\in\mathbb Z_+.$ If $R_1$ is integral over $R_2$ then
$$\lim_{m\to\infty}\lambda(I(1)_{mn_1}\cdots I(r)_{mn_r}/J(1)_{mn_1}\cdots J(r)_{mn_r})/m^d=0.$$
Converse holds if $R$ is formally equidimensional and $\grade(\J)\geq 1.$
\et
\begin{proof}
	Since for all $1\leq l\leq r,$ $\mathcal J(l)=\{J(l)_n\}$ are Noetherian filtrations of ideals in $R,$ there exists an integer $\alpha\in\mathbb Z_+,$ such that $G_l^{\alpha}=\oplus_{n\geq 0}J(l)_{\alpha n}$ are Noetherian standard graded $R$-algebras for all $1\leq l\leq r.$ Hence $T=\bigoplus_{m_1,\ldots,m_r\geq 0}J(1)_{\alpha m_1}\cdots J(r)_{\alpha m_r}$ is a Noetherian standard graded $R$-algebra and graded subring of $R_2$. 
	\\Using Lemma \ref{use} for $u_1=\ldots=u_r=\alpha,$ we get  $R_2$ is finitely generated $T$-module. Since $R_1$ is finitely generated $R_2$-module, we have $R_1$ is finitely generated $T$-module. Now for all integers $0\leq b_j\leq \alpha -1$ with $1\leq j\leq r,$ $S^{(\alpha,b_1,\ldots,b_r)}=\bigoplus_{m_1,\ldots,m_r\geq 0}I(1)_{\alpha m_1+b_1}\cdots I(r)_{\alpha m_r+b_r}$ are $T$-submodules of $R_1.$ Therefore  $S^{(\alpha,b_1,\ldots,b_r)}$ are finitely generated $T$-modules for all  integers $0\leq b_j\leq \alpha -1$ with $1\leq j\leq r$ and hence $$G^{(b_1,\ldots,b_r)}=\bigoplus_{m_1,\ldots,m_r\geq 0}I(1)_{\alpha m_1+b_1}\cdots I(r)_{\alpha m_r+b_r}/J(1)_{\alpha m_1+b_1}\cdots J(r)_{\alpha m_r+b_r}$$ are finitely generated $T$-modules for all  integers $0\leq b_j\leq \alpha -1$ with $1\leq j\leq r$ where we consider the grading $$G^{(b_1,\ldots,b_r)}_{m_1,\ldots,m_r}=I(1)_{\alpha m_1+b_1}\cdots I(r)_{\alpha m_r+b_r}/J(1)_{\alpha m_1+b_1}\cdots J(r)_{\alpha m_r+b_r}$$ and $T_{m_1,\ldots,m_r}=J(1)_{\alpha m_1}\cdots J(r)_{\alpha m_r}.$ Since  $\lambda( I(l)_n/J(l)_n)<\infty$ for all $n\in\mathbb N$ and $1\leq l\leq r,$ there exists an integer $t\in\mathbb Z_+,$ such that $G^{(b_1,\ldots,b_r)}$ is finitely generated $T/\mm^tT$-module for all $0\leq b_1,\ldots,b_r\leq \alpha-1.$
	\\ By  \cite[Theorem 4.1]{HHRT}, for all $1\leq j\leq r$ and integers $0\leq b_j\leq \alpha-1,$ there exist polynomials $P_{(b_1,\ldots,b_r)}(X_1,\ldots,X_r)\in\mathbb Q[X_1,\ldots,X_r]$  of total degree at most $\dim \supp_{++}T/\mm^tT$ and an integer $f\in\mathbb Z_+$ such that for all $m_1,\ldots,m_r\geq f,$ $$P_{(b_1,\ldots,b_r)}(m_1,\ldots,m_r)=\lambda(G^{(b_1,\ldots,b_r)}_{m_1,\ldots,m_r}).$$
	\\Since $J(1)_{\alpha m_1}\cdots J(r)_{\alpha m_r}=J(1)_{\alpha }^{m_1}\cdots J(r)_{\alpha}^{ m_r},$ by \cite{KS}, \cite[Corollary 3.3 and Corollary 5.3]{F}, we get 
	\begin{equation}\label{dimension}  \dim \supp_{++}T/\mm^tT\leq\dim \supp_{++}T/\mm T=l(J(1)_{\alpha }\cdots J(r)_{\alpha})-1<d\end{equation}
	where $l(J(1)_{\alpha }\cdots J(r)_{\alpha})$ is the analytic spread of $J(1)_{\alpha }\cdots J(r)_{\alpha}.$
	\\Let $i_1,\ldots,i_r\in \mathbb N$ with $i_1+\cdots+i_r< d$ and $$P_{(b_1,\ldots,b_r)}(X_1,\ldots,X_r)=\sum_{i_1+\cdots+i_r< d}z_{i_1,\ldots,i_r}(b_1,\ldots,b_r)X_1^{i_1}X_2^{i_2}\cdots X_r^{i_r}$$ where $z_{i_1,\ldots,i_r}(b_1,\ldots,b_r)\in\mathbb Q.$ Let {\small\begin{equation*}\begin{array}{lll}&&C:=\max\{|z_{i_1,\ldots,i_r}(b_1,\ldots,b_r)|(n_1\cdots n_r)^d, \lambda(I(1)_{mn_1}\cdots I(r)_{mn_r}/J(1)_{mn_1}\cdots J(r)_{mn_r}) \\&&~~~~~~~~~~~~~~:0\leq b_1,\ldots,b_r\leq \alpha-1, \sum\limits_{j=1}^ri_j< d \mbox{ and } 0\leq m\leq \alpha (f+1) \}.\end{array}\end{equation*}}
	Then for any $m\in\mathbb Z_+,$ by equation (\ref{dimension}), we get  $$\lambda(I(1)_{mn_1}\cdots I(r)_{mn_r}/J(1)_{mn_1}\cdots J(r)_{mn_r})/m^d\leq C/m<C.$$ Therefore 
	\begin{equation*}\label{eq}
	\begin{array}{lll}
	0&\leq& \liminf_{m\to\infty}\lambda(I(1)_{mn_1}\cdots I(r)_{mn_r}/J(1)_{mn_1}\cdots J(r)_{mn_r})/m^d\\&\leq&\limsup_{m\to\infty}\lambda(I(1)_{mn_1}\cdots I(r)_{mn_r}/J(1)_{mn_1}\cdots J(r)_{mn_r})/m^d=0\end{array}\end{equation*} implies $$\lim_{m\to\infty}\lambda(I(1)_{mn_1}\cdots I(r)_{mn_r}/J(1)_{mn_1}\cdots J(r)_{mn_r})/m^d=0.$$
	\\Now we prove the converse. Let $R$ be a formally equidimensional local ring, \\$\grade(\J)\geq 1$ and $\lim_{m\to\infty}\lambda(I(1)_{mn_1}\cdots I(r)_{mn_r}/J(1)_{mn_1}\cdots J(r)_{mn_r})/m^d=0.$ 
	\\Suppose $R_1$ is not integral over $R_2.$
	\\{\bf{Claim :}} For any  $u_1,\ldots,u_r\in\mathbb Z_+,$ $R_1^{(u_1,\ldots,u_r)}=\oplus_{m_1,\ldots,m_r\geq 0}I(1)_{u_1m_1}\cdots I(r)_{u_rm_r}$ is not integral over $R_2^{(u_1,\ldots,u_r)}=\oplus_{m_1,\ldots,m_r\geq 0}J(1)_{u_1m_1}\cdots J(r)_{u_rm_r}.$
	\\{\bf{Proof of the claim:}} Suppose there exist  $u_1,\ldots,u_r\in\mathbb Z_+$ such that $R_1^{(u_1,\ldots,u_r)}$ is integral over $R_2^{( u_1,\ldots,u_r)}.$ By Lemma \ref{use}, $R_1$ is finitely generated $R_1^{(u_1,\ldots,u_r)}$-module. Hence $R_1$ is integral over  $R_2^{(u_1,\ldots,u_r)}$ as well as over $R_2$ which is a contradiction.
	\\Since for all $1\leq l\leq r,$ $\mathcal I(l)=\{I(l)_n\},$ $\mathcal J(l)=\{J(l)_n\}$ are Noetherian filtrations of ideals in $R,$ there exists an integer $\alpha\in\mathbb Z_+,$ such that $\oplus_{n\geq 0}I(l)_{\alpha n},$ $\oplus_{n\geq 0}J(l)_{\alpha n}$ are Noetherian standard graded $R$-algebras for all $1\leq l\leq r.$ Now by the claim, 
	$$R_1^{(\alpha n_1,\ldots,\alpha n_r)}=\oplus_{m_1,\ldots,m_r\geq 0}I(1)_{\alpha n_1m_1}\cdots I(r)_{\alpha n_rm_r}=\oplus_{m_1,\ldots,m_r\geq 0}I(1)_{\alpha n_1}^{m_1}\cdots I(r)_{\alpha n_r}^{m_r}$$
	is not integral over 
	$$R_2^{(\alpha n_1,\ldots,\alpha n_r)}=\oplus_{m_1,\ldots,m_r\geq 0}J(1)_{\alpha n_1m_1}\cdots J(r)_{\alpha n_rm_r}=\oplus_{m_1,\ldots,m_r\geq 0}J(1)_{\alpha n_1}^{m_1}\cdots J(r)_{\alpha n_r}^{m_r}.$$
	Therefore there exists $t\in\{1,\ldots,r\}$ such that $\oplus_{m_{t}\geq 0}I(t)_{\alpha n_{t}}^{m_{t}}$ is not finitely generated over $\oplus_{m_{t}\geq 0}J(t)_{\alpha n_{t}}^{m_{t}},$  i.e., $J(t)_{\alpha n_{t}}$ is not a reduction of $I(t)_{\alpha n_{t}}.$ 
	\\Let $n=\max\{n_1,\ldots,n_r\}.$ Since $(\J)^{n\alpha}\subset J(1)_{\alpha n_1}\cdots J(r)_{\alpha n_r},$ by Lemma \ref{reduction}, $J(1)_{\alpha n_1}\cdots J(r)_{\alpha n_r}$ is not a reduction of $I(1)_{\alpha n_1}\cdots I(r)_{\alpha n_r}.$ Therefore 
	$$\mathcal R(I(1)_{\alpha n_1}\cdots I(r)_{\alpha n_r})=\oplus_{m\geq 0}(I(1)_{\alpha n_1}\cdots I(r)_{\alpha n_r})^{m}=\oplus_{m\geq 0}I(1)_{\alpha n_1m}\cdots I(r)_{\alpha n_rm}$$ is not finitely generated module over
	$$\mathcal R(J(1)_{\alpha n_1}\cdots J(r)_{\alpha n_r})=\oplus_{m\geq 0}(J(1)_{\alpha n_1}\cdots J(r)_{\alpha n_r})^{m}=\oplus_{m\geq 0}J(1)_{\alpha n_1m}\cdots J(r)_{\alpha n_rm}.$$
	Thus by \cite[Theorem 2.1]{Rees}, $\lim_{m\to\infty}\lambda(I(1)_{\alpha mn_1}\cdots I(r)_{\alpha mn_r}/J(1)_{\alpha mn_1}\cdots J(r)_{\alpha mn_r})/m^d\neq 0,$ which is a contradiction.
	\end{proof}	
Now we prove the ``only if" part of Rees' theorem \ref{Reestheorem} for (not necessarily Noetherian) filtrations. Let $(R,\mm)$ be a complete local domain of dimension $d>0.$ Then by \cite[Lemma 4.2]{CSS},  \cite{C2014}, there exists a regular local ring $S$ of dimension $d$ which birationally dominates $R$. Let $\lm_1,\ldots,\lm_d\in \mathbb R$ be rationally independent real numbers such that $\lambda_i\ge 1$ for all $i$ and  $y_1,\ldots,y_d\in S$ be a regular system of parameters in $S.$ We consider a valuation $\nu$ on the quotient field of $R$ which dominates $S$ as mentioned in \cite{C2014}, i.e., $\nu(y_1^{a_1}\cdots y_d^{a_d})=a_1\lambda_1+\cdots+a_d\lambda_d$ for $a_1,\ldots,a_d\in \NN$ and $\nu(\gamma)=0$ if $\gamma\in S$ is a unit. Let $k=R/m_R$ and $k'=S/m_S$. Then by \cite{C2014} (Page 10), we get $k'=S/\mm_S=V_{\nu}/\mm_{\nu}$ where $V_{\nu}$ is the valuation ring of $\nu$ and $[k':k]<\infty.$

For $\mu\in \mathbb R_{>0}$, define ideals $K_{\mu}$ and $K_{\mu}^+$ in the valuation ring $V_{\nu}$  by
$$
K_{\mu}=\{f\in {\rm Q}(R)\mid \nu(f)\ge \mu\}, 
$$
$$
K_{\mu}^+=\{f\in {\rm Q}(R)\mid \nu(f)>\mu\}.
$$
Let $\mathcal T(m)=\{T(m)_n\}$ be filtrations of ideals in $R$ for all $1\leq m\leq r$ and $a\in\mathbb Z_+.$ We recall the definition of $a$th truncated filtration of ideals defined in \cite{CSS}. The $a$-{\it {th truncated filtration}} of ideals  $\mathcal T_a(m)=\{T_a(m)_n\}$ of $\mathcal T(m)$ is defined  by $T_a(m)_n=T(m)_n$ if $n\le a$ and if $n>a$, then $T_a(m)_n=\sum T_a(m)_iT_a(m)_j$ where the sum is over $i,j>0$ such that $i+j=n$. 

 Fix $n_1,\ldots,n_r\in\mathbb Z_+.$ Define filtrations  of ideals 
$$\mathcal T=\{\mathcal T_m:=T(1)_{mn_1}\cdots T(r)_{mn_r}\}\mbox{ and }\mathcal T(a)=\{\mathcal T(a)_m:=T_a(1)_{mn_1}\cdots T_a(r)_{mn_r}\}.$$ By \cite[Lemma 4.3]{C2013}, there exists $\alpha\in \ZZ_+$ such that $K_{\alpha m}\cap R\subset \mm^m$ for all $m\in \NN.$ 
\\ Note that $\mm^{\alpha bmd}\mathcal T_m\subset \mathcal T_m\cap K_{\alpha mb}$ and $\mm^{\alpha bmd}\mathcal T(a)_m\subset \mathcal T(a)_m\cap K_{\alpha mb}$ for all $m\in \NN$ and $b\in\mathbb Z_+.$ 
\bp\label{Noetherianneed}
	Suppose $(R,\mm),$ $\mathcal T,$ $\mathcal T(a),$ $\alpha$ and $K_{\mu}$ as above. Fix $b'\in\mathbb Z_+.$ For each $a\in\mathbb Z_+,$ let $\mathcal F(a)=\{  \mathcal F(a)_n\}$ be a filtration of ideals such that for all $n,$ $$\mathcal T(a)_n\subset  \mathcal F(a)_n\subset \mathcal T_n.$$	Then there exists an integer $b$ such that $b\geq b'$ and \begin{equation*}\begin{array}{lll}&&\lim_{a\to\infty}\big({\lim_{m\to\infty}\lambda({\mathcal T(a)_m/\mathcal T(a)_m\cap K_{\alpha mb}})}\big)/m^d=\lim_{m\to\infty}\lambda({\mathcal T_m/\mathcal T_m\cap K_{\alpha mb})/m^d}\\&=&\lim_{a\to\infty}\big({\lim_{m\to\infty}\lambda({\mathcal F(a)_m/\mathcal F(a)_m\cap K_{\alpha mb}})}\big)/m^d.\end{array}\end{equation*}	
\ep 
\begin{proof}
	We follow the argument given in \cite{C2014}, \cite[Proposition 4.3]{CSS}. 
	For $t\geq 1$, define the semigroups
	$$
	\begin{array}{lll}
	\Gamma'^{(t)}&=&\{(m_1,\ldots,m_d,i)\in \NN^{d+1}\mid \dim_k  \mathcal T_i\cap K_{m_1\lambda_1+\cdots+m_d\lambda_d}/\mathcal T_i\cap K^+_{m_1\lambda_1+\cdots+m_d\lambda_d} \ge t\\
	&&\mbox{ and }m_1+\cdots+m_d\le \alpha b'i\},
	\end{array}
	$$
	$$
	\begin{array}{lll}
	\Gamma(a)'^{(t)}&=&
	\{(m_1,\ldots,m_d,i)\in \NN^{d+1}\mid \\ &&\dim_k \mathcal T(a)_i\cap K_{m_1\lambda_1+\cdots+m_d\lambda_d}/\mathcal T(a)_i\cap K^+_{m_1\lambda_1+\cdots+m_d\lambda_d} \ge t\\
	&&\mbox{ and }m_1+\cdots+m_d\le \alpha b'i\}
	\end{array}
	$$ 
	$$
	\begin{array}{lll}
	\Gamma(\mathcal F(a))'^{(t)}&=&
	\{(m_1,\ldots,m_d,i)\in \NN^{d+1}\mid \\&&\dim_k \mathcal F(a)_i\cap K_{m_1\lambda_1+\cdots+m_d\lambda_d}/\mathcal F(a)_i\cap K^+_{m_1\lambda_1+\cdots+m_d\lambda_d} \ge t\\
	&&\mbox{ and }m_1+\cdots+m_d\le \alpha b'i\}.
	\end{array}
	$$ 
	By  \cite[Lemma 4.5]{C2014}, there exists an integer $b\geq b'$ such that the semigroups 
	$$
	\begin{array}{lll}
	\Gamma^{(t)}&=&\{(m_1,\ldots,m_d,i)\in \NN^{d+1}\mid \dim_k  \mathcal T_i\cap K_{m_1\lambda_1+\cdots+m_d\lambda_d}/\mathcal T_i\cap K^+_{m_1\lambda_1+\cdots+m_d\lambda_d} \ge t\\
	&&\mbox{ and }m_1+\cdots+m_d\le \alpha bi\},
	\end{array}
	$$
	$$
	\begin{array}{lll}
	\Gamma(a)^{(t)}&=&
	\{(m_1,\ldots,m_d,i)\in \NN^{d+1}\mid \\ &&\dim_k \mathcal T(a)_i\cap K_{m_1\lambda_1+\cdots+m_d\lambda_d}/\mathcal T(a)_i\cap K^+_{m_1\lambda_1+\cdots+m_d\lambda_d} \ge t\\
	&&\mbox{ and }m_1+\cdots+m_d\le \alpha bi\}
	\end{array}
	$$ 
	$$
	\begin{array}{lll}
	\Gamma(\mathcal F(a))^{(t)}&=&
	\{(m_1,\ldots,m_d,i)\in \NN^{d+1}\mid \\&&\dim_k \mathcal F(a)_i\cap K_{m_1\lambda_1+\cdots+m_d\lambda_d}/\mathcal F(a)_i\cap K^+_{m_1\lambda_1+\cdots+m_d\lambda_d} \ge t\\
	&&\mbox{ and }m_1+\cdots+m_d\le \alpha bi\}.
	\end{array}
	$$ satisfy equations $(5)$ and $(6)$ of \cite{C2014}.
	
	Define $\Gamma_m^{(t)}=\Gamma^{(t)}\cap (\NN^d\times\{m\}),$ $\Gamma(a)^{(t)}_m=\Gamma(a)^{(t)}\cap (\NN^d\times\{m\})$ and $\Gamma(\mathcal F(a))^{(t)}_m=\Gamma(\mathcal F(a))^{(t)}\cap (\NN^d\times\{m\})$  for $m\in \NN$. For a (strongly nonnegative) sub semigroup $S$ of $\ZZ^d\times \NN,$ ${\rm con}(S)$ is defined as the closed convex cone which is the closure of the set of all linear combinations $\sum\lambda_i s_i$  with $s_i\in S$ and $\lambda_i$ a nonnegative real number and  the {\it{Newton-Okounkov body}}  is defined as
	$$
	\Delta(S)={\rm con}(S)\cap (\mathbb R^d\times \{1\}).
	$$
	This theory is developed in \cite{Ok}, \cite{LM} and \cite{KK} and is summarized in  \cite[Section 3]{C2014}. 
By \cite[Lemma 4.5]{C2014} and \cite[Theorem 3.2]{C2014}, 
	\begin{equation}\label{eq3}
	\lim_{m\rightarrow\infty}\frac{\#\Gamma_m^{(t)}}{m^d}={\rm Vol}(\Delta(\Gamma^{(t)})),
	\end{equation}
	\begin{equation}\label{eq4}
	\lim_{m\rightarrow\infty}\frac{\#\Gamma(a)_m^{(t)}}{m^d}={\rm Vol}(\Delta(\Gamma(a)^{(t)})),
	\end{equation}
	\begin{equation}\label{eql4}
	\lim_{m\rightarrow\infty}\frac{\#\Gamma(\mathcal F(a))_m^{(t)}}{m^d}={\rm Vol}(\Delta(\Gamma(\mathcal F(a))^{(t)})).
	\end{equation}
	Therefore
	\begin{equation}\label{eqb1}
	\begin{array}{lll}
	\lim_{m\to\infty}\lambda({\mathcal T_m/\mathcal T_m\cap K_{\alpha mb})/m^d}&=& \lim_{m\rightarrow\infty}\dim_k\big(\oplus_{0\leq\mu<\alpha mb}\mathcal T_m\cap K_{\mu}/\mathcal T_m\cap K_{\mu}^+\big)/m^d\\&=&\sum_{t=1}^{[k':k]}\lim_{m\rightarrow\infty}\frac{\#\Gamma^{(t)}_m}{m^d}.\end{array}\end{equation}
	Similarly
	\begin{equation}\label{eq2}
\begin{array}{lll}
	&&\lim_{m\to\infty}\lambda({\mathcal T(a)_m/\mathcal T(a)_m\cap K_{\alpha mb})/m^d}\\&=& \lim_{m\rightarrow\infty}\dim_k\big(\oplus_{0\leq\mu<\alpha mb}\mathcal T(a)_m\cap K_{\mu}/\mathcal T(a)_m\cap K_{\mu}^+\big)/m^d
	\\&=&\sum_{t=1}^{[k':k]}\lim_{m\rightarrow\infty}\frac{\#\Gamma(a)^{(t)}_m}{m^d}
	\end{array}
	\end{equation}
	and
	\begin{equation}\label{eq2}
	\begin{array}{lll}
	&&\lim_{m\to\infty}\lambda({\mathcal F(a)_m/\mathcal F(a)_m\cap K_{\alpha mb})/m^d}\\&=&\lim_{m\rightarrow\infty}\dim_k\big(\oplus_{0\leq\mu<\alpha mb}\mathcal F(a)_m\cap K_{\mu}/\mathcal F(a)_m\cap K_{\mu}^+\big)/m^d\\&=&\sum_{t=1}^{[k':k]}\lim_{m\rightarrow\infty}\frac{\#\Gamma(\mathcal F(a))^{(t)}_m}{m^d}.\end{array}\end{equation}
	Let $\hat a=\lfloor a/\max\{n_1,\ldots,n_r\}\rfloor$ where $\lfloor x \rfloor$ be the largest integer smaller than or equal to the real number $x.$ Now \begin{equation}\label{eq5}
	\Gamma_i^{(t)}=\Gamma(\mathcal F(a))_i^{(t)}=\Gamma(a)_i^{(t)}\mbox{ for }i\le \hat a.
	\end{equation}
	Hence
	$$
	n*\Gamma_{\hat a}^{(t)}:=\{x_1+\cdots+x_n\mid x_1,\ldots,x_n\in \Gamma_{\hat a}^{(t)}\} \subset \Gamma(a)^{(t)}_{n\hat a}\subset \Gamma(\mathcal F(a))^{(t)}_{n\hat a}\mbox{ for all }n\ge 1.
	$$
	By  \cite[Theorem 3.3]{C2014} and since $\hat a\to\infty$ as $a\to\infty,$ given $\epsilon>0$, there exists $a_0>0$ such that for all $a\ge a_0,$ we have
	\begin{equation}\label{eq6}
	\begin{array}{lll}
	{\rm Vol}(\Delta(\Gamma^{(t)}))&\ge &{\rm Vol}(\Delta(\Gamma(a)^{(t)})
	=\lim_{n\rightarrow\infty}\frac{\#\Gamma(a)^{(t)}_n}{n^d} =\lim_{n\rightarrow\infty}\frac{\#\Gamma(a)^{(t)}_{n\hat a}}{(n\hat a)^d}\\
	&\ge &\lim_{n\rightarrow\infty}\frac{\#(n*\Gamma_{\hat a}^{(t)})}{(n\hat a)^d}\ge {\rm Vol}(\Delta(\Gamma^{(t)}))-\epsilon
	\end{array}
	\end{equation}
	and
	\begin{equation}\label{eq6}
	\begin{array}{lll}
	{\rm Vol}(\Delta(\Gamma^{(t)}))&\ge &{\rm Vol}(\Delta(\Gamma(\mathcal F(a))^{(t)})
	=\lim_{n\rightarrow\infty}\frac{\#\Gamma(\mathcal F(a))^{(t)}_n}{n^d} =\lim_{n\rightarrow\infty}\frac{\#\Gamma(\mathcal F(a))^{(t)}_{n\hat a}}{(n\hat a)^d}\\
	&\ge &\lim_{n\rightarrow\infty}\frac{\#(n*\Gamma_{\hat a}^{(t)})}{(n\hat a)^d}\ge {\rm Vol}(\Delta(\Gamma^{(t)}))-\epsilon.
	\end{array}
	\end{equation}
	By equations (\ref{eqb1})-(\ref{eq6}), we get the desired result.
\end{proof}
\bt\label{domain}
	Let $(R,\mm)$ be an analytically irreducible local ring of dimension $d>0$ and for $1\leq l\leq r,$ $\mathcal I(l)=\{I(l)_n\},$ $\mathcal J(l)=\{J(l)_n\}$ be (not necessarily Noetherian) filtrations of ideals in $R$ such that $J(l)_n\subset I(l)_n$ for all $n\in\mathbb N$ and $1\leq l\leq r.$ Fix $n_1,\ldots,n_r\in\mathbb Z_+.$ Suppose there exists an integer $c\in\mathbb Z_+$ such that for all $i\in\mathbb N,$$$\mm^{ci}\cap I(1)_{in_1}\cdots I(r)_{in_r}=\mm^{ci}\cap J(1)_{in_1}\cdots J(r)_{in_r}. $$ Let $R_1=\oplus_{m_1,\ldots,m_r\geq 0}I(1)_{m_1}\cdots I(r)_{m_r}$ be integral over $R_2=\oplus_{m_1,\ldots,m_r\geq 0}J(1)_{m_1}\cdots J(r)_{m_r}.$ Then $$\lim_{m\to\infty}\lambda(I(1)_{mn_1}\cdots I(r)_{mn_r}/J(1)_{mn_1}\cdots J(r)_{mn_r})/m^d=0.$$
\et
\begin{proof}
	Let $\hat R$ denote the $\mm$-adic completion of $R.$ Without loss of generality, we may replace $R$ by $\hat R$ and $I(l)_n, J(l)_n$ by $I(l)_n{\hat R},J(l)_n{\hat R}$ respectively for all $n$ and  $1\leq l\leq r.$
	\\For $a\in\mathbb Z_+,$ we define filtrations of ideals 
	$$\mathcal I=\{\mathcal I_i:=I(1)_{in_1}\cdots I(r)_{in_r}\}, \mbox{ } \mathcal I(a)=\{\mathcal I(a)_i:=I_a(1)_{in_1}\cdots I_a(r)_{in_r}\},$$ and $$\mathcal J=\{\mathcal J_i:=J(1)_{in_1}\cdots J(r)_{in_r}\}.$$  Now there exists an integer $\underline a\in\mathbb Z_+,$	such that $a\leq \underline a$ and every element of $\oplus_{i\geq 0}\mathcal I(a)_i$ is integral over $\oplus_{i\geq 0}\mathcal J(\underline a)_i$ where $$\mathcal J(\underline a)=\{\mathcal J(\underline  a)_i:=J_{\underline a}(1)_{in_1}\cdots J_{\underline a}(r)_{in_r}\}.$$
	Define the following filtrations of ideals
	$$\mathcal A_{a}=\{\mathcal A_{a,m}:=\sum_{\substack{\alpha+\beta=m\\ \alpha,\beta\in\mathbb N}}\mathcal I(a)_\alpha{\mathcal  J(\underline a)_\beta\}}$$ and $$\mathcal B_{\underline a}=\{\mathcal B_{\underline a,m}:=\mathcal A_{a,m}\cap \mathcal J_m\}.$$ Note that $\oplus_{n\geq 0}\mathcal A_{a,n}$ is finitely generated $\oplus_{n\geq 0}\mathcal B_{\underline a,n}$-module. Thus by Theorem \ref{Noetheriancase}, \begin{equation}\label{eq7}\lim_{m\to\infty}\lambda(\mathcal A_{a,m}/\mathcal B_{\underline a,m})/m^d=0.\end{equation}
	Now for all $n,$
	$$	\mathcal A_{a,n}\cap \mm^{cn}=\mathcal A_{a,n}\cap\mathcal I_n\cap \mm^{cn}=\mathcal A_{a,n}\cap \mathcal J_n\cap \mm^{cn}=\mathcal B_{\underline a,n}\cap \mm^{cn}.$$
 Let $\alpha\in \ZZ_+$ be the integer such that $K_{\alpha n}\cap R\subset \mm^n$ for all $n\in \NN$ \cite[Lemma 4.3]{C2013}. Therefore 
	$$\mathcal A_{a,n}\cap K_{\alpha cn}=\mathcal A_{a,n}\cap K_{\alpha cn}\cap R\subset \mathcal A_{a,n}\cap \mm^{cn}\subset\mathcal B_{\underline a,n}$$ and hence $\mathcal A_{a,n}\cap K_{\alpha cn}=\mathcal B_{\underline a,n}\cap K_{\alpha cn}$	for all $n.$
	
	Now for all $n\in\mathbb N,$ $$\mathcal I(a)_n\subset \mathcal A_{a,n}\subset  \mathcal I_n, \mbox{  }\mathcal J(\underline a)_n\subset \mathcal B_{\underline a,n}\subset  \mathcal J_n.$$ 
	Choose an integer $c'\geq c$ such that Proposition \ref{Noetherianneed} holds for the above two equations. Also note that $i\in\mathbb N,$$$\mm^{c'i}\cap I(1)_{in_1}\cdots I(r)_{in_r}=\mm^{c'i}\cap J(1)_{in_1}\cdots J(r)_{in_r}. $$ Replace $c$ by $c'$. Since $\underline a\to\infty$ as $a\to\infty,$ by equation (\ref{eq7}) and Proposition \ref{Noetherianneed},  we get 
	\begin{equation*}
	\begin{array}{lll}	
	&&\lim_{m\to\infty}\lambda(\mathcal I_m/\mathcal J_m)/m^d\\&=&\lim_{m\to\infty}\lambda(\mathcal I_m/\mathcal I_m\cap K_{\alpha cm})/m^d-\lim_{m\to\infty}\lambda(\mathcal J_m/\mathcal J_m\cap K_{\alpha cm})/m^d\\&=&\lim_{a\to\infty}\lim_{m\to\infty}\lambda(\mathcal  A_{a,m}/\mathcal A_{a,m}\cap K_{\alpha cm})/m^d\\&&-\lim_{a\to\infty}\lim_{m\to\infty}\lambda(\mathcal B_{\underline a,m}/\mathcal B_{\underline a,m}\cap K_{\alpha cm})/m^d\\&=&\lim_{a\to\infty}\lim_{m\to\infty}\lambda(\mathcal A_{a,m}/\mathcal B_{\underline a,m})/m^d=0.
	\end{array}
	\end{equation*}
\end{proof}
It is natural to ask the following.
\begin{question}
		Does Theorem \ref{domain} hold true for any analytically unramified local ring $?$
	\end{question}
The author believes the answer is positive but is unable to prove.
\\The following example shows that the ``if" part of Rees' theorem \ref{Reestheorem} $(2)$ is not true for non-Noetherian filtrations of ideals in general.
\bex\label{counter}
	Let $R=k[|x,y|]$ be a formal power series ring of dimension two over a field $k.$ Let $\mathcal I=\{I_m:=x^m\}$ and $\mathcal J=\{J_m:=(x^{m+1},x^my)\}.$ Note that $I_m\cap \mm^{2m}=J_m\cap \mm^{2m}$ for all $m\geq 0$ and $\lim_{m\to\infty}\lambda( I_{mn}/J_{mn})/m^2=0$ but $\oplus_{m\geq 0}I_m$ is not integral over $\oplus_{m\geq 0}J_m.$
\eex
Let $r\geq 2.$ We define a {\it{multigraded filtration}} $\mathcal I=\{I_{\n}\}$  where $\n=(n_1,\ldots,n_r)\in\NN^r$ of ideals in a ring $R$ to be a collection of ideals of $R$ such that
$R=I_{\bf 0}$, $I_{\n+{\bf e_i}}\subset I_{\n}$
for all $i=1,\ldots,r$ where $e_i=(0,\ldots,1,\ldots,0)$ with $1$ at $i$th position and
$I_{\n}I_{\m}\subset I_{\n+\m}$ whenever $\n,\m\in \NN^r$. As a generalization of Rees' theorem \ref{Reestheorem}, we pose the following question for  multigraded filtration of ideals.
\begin{question}
	Let  $\{I_{\n}\}$ and  $\{J_{\n}\}$ be multigraded filtrations of ideals in a Noetherian local ring such that $J_{\n}\subset I_{\n}$ and $\lm( I_{\n}/ J_{\n})<\infty$ for all $\n\in\NN^r.$ Fix $n_1,\ldots,n_r\in\ZZ_+.$
	\ben
	\item  Is it true that if  $\bigoplus\limits_{\n\in\NN^r}I_{\n}$ is integral over $\bigoplus\limits_{\n\in\NN^r}J_{\n}$ then $$\lim\limits_{m\to\infty}\lambda(I_{mn_1,\ldots,mn_r}/J_{mn_1,\ldots,mn_r})/m^d=0 ?$$
	\item Does the converse of $(1)$ hold true if $\bigoplus\limits_{\n\in\NN^r}I_{\n}$ and $\bigoplus\limits_{\n\in\NN^r}J_{\n}$ both are Noetherian $?$ 
	\een
	\end{question}

\section{the degree of the rees polynomial}
In this section, we provide bounds for the degree of the function $\lm(I^n/J^n)$ where $J\subsetneq I$ are ideals in a Noetherian local ring with $\lm(I/J)<\infty$ and $J$ is reduction of $I.$ We concentrate on the ideals generated by $d$-sequences. The notion of $d$-sequence was introduced by Huneke \cite{Hu}. Let ${\underline x}=x_1,\ldots,x_n$ be a sequence of elements in $R.$ Then ${\underline x}$ is called a {\it{$d$-sequence}} if
\ben
\item[(1)] $x_i\notin (x_1,\ldots,{\widehat{x_i}},\ldots,x_n)$ for all $i=1,\ldots,n,$
\item[(2)]
for all $k\geq i+1$ and all $i\geq 0,$ $(x_0=0)$
$$((x_0,\ldots,x_i):x_{i+1}x_k)=((x_0,\ldots,x_i):x_k).$$
\een Every regular sequence is a $d$-sequence but the converse is not true. Every system of parameters in a Buchsbaum local ring is a $d$-sequence. Examples of $d$-sequences are abundant. Ideals generated by $d$-sequences have nice properties, e.g. they are of linear type \cite{Hu}, Castelnuovo-Mumford regularity of Rees algebras of such ideals is zero \cite{Trung}. For our convenience we omit the condition $(1)$ in the definition of $d$-sequence. We say, a sequence of elements ${\underline x}=x_1,\ldots,x_n$ in $R$ is a {\it{$d$-sequence}} if condition $(2)$ is satisfied.
\begin{remark}\label{n2}
	Let $J\subsetneq I$ be ideals in a Noetherian local ring $(R,\mm)$ with the property that $\lm(I/J)<\infty.$ Suppose $J$ is a reduction of $I$ and $JI^m=I^{m+1}$ for all $m\geq l.$ Consider the Veronese subring $A=\bigoplus_{n\geq 0}J^{nl}$ of $\R(J).$ Then $M=\bigoplus_{n\geq 0}I^{nl}/J^{nl}$ is a finite $A$-module. Therefore there exists an integer $t\in\ZZ_+,$ such that 
$\mm^tI^{nl}\subseteq J^{nl}$ for all $n\geq 1.$ Hence $M$ is $A/\mm^tA$-module. Thus $\dim M\leq\dim A/\mm A=\dim \R(J)/\mm\R(J)=l(J).$ Therefore $\lm (I^{nl}/J^{nl})$ is a polynomial type function of degree at most $l(J)-1.$ Hence $\lm(I^n/J^n)$ is a polynomial type function of degree at most $l(J)-1.$
\end{remark} 
Let $J\subset I$ be ideals in a Noetherian local ring $(R,\mm)$ and $a_1,\ldots,a_s$ be a sequence of elements in $J.$ Throughout this section we use the following notation: $R_i=R/(a_0,a_1,\ldots,a_i)$ for all $i=0,\ldots,s$ where $a_0=0$ and $I_i,J_i$ denote the images of $I,J$ in $R_i$ respectively for all $i=0,\ldots,s.$ The next two results are requied to prove our main results.
\bl\label{need}
Let $(R,\mm)$ be a Noetherian local ring, $J\subsetneq I$ be ideals in $R$ and $\lm(I/J)<\infty.$ Let $a_1,\ldots,a_s$ be a sequence of elements in $J.$ Suppose for some $i\in \{0,1,\ldots,s-1\}$ and all $n\gg 0,$ $ (J^{n+1}_i:(a^{(i)}_{i+1}))\cap I^{n}_i=J^{n}_i$ where  $a^{(i)}_{i+1}$ denotes the image of $a_{i+1}$ in $R_i.$ Then for all $n\gg0,$ $$\lm(I^{n+1}_{i+1}/J^{n+1}_{i+1})\leq \lm(I^{n+1}_i/J^{n+1}_i)-\lm(I^{n}_i/J^{n}_i) .$$ \el
\bpf
For all $n\gg 0,$ consider the following exact sequence {\small{$$0\longrightarrow(J^{n+1}_i:(a^{(i)}_{i+1}))\cap I^{n}_i/J^{n}_i\longrightarrow I^{n}_i/J^{n}_i\overset{.a^{(i)}_{i+1}}\longrightarrow I^{n+1}_i/J^{n+1}_i\longrightarrow I^{n+1}_i/a^{(i)}_{i+1}I^{n}_i+J^{n}_i\longrightarrow 0.$$} } 
Since $$I^{n+1}_{i+1}/J^{n+1}_{i+1}={I^{n+1}_{i}+(a^{(i)}_{i+1})}/{J^{n+1}_{i}+(a^{(i)}_{i+1})}\simeq I^{n+1}_{i}/(a^{(i)}_{i+1})\cap I^{n+1}_{i}+J^{n+1}_{i}$$ and $a^{(i)}_{i+1}I^{n}_{i}+J^{n+1}_{i}\subseteq (a^{(i)}_{i+1})\cap I^{n+1}_{i}+J^{n+1}_{i},$ we have the required inequality.
\eepf
\bp\label{n1}
Let $x_1,\ldots,x_n$ be a $d$-sequence in a Noetherian local ring $(R,\mm)$ and $J=(x_1,\ldots,x_n).$ Then $J^n\cap (x_1)=x_1J^{n-1}$ for all $n\geq 1.$
\ep
\bpf
Note that images of $x_2,\ldots,x_n$ is a $d$-sequence in $R/(x_1).$ Let $x_1r\in J^n\cap (x_1).$ Then $x_1r=x_1a+b$ for some $a\in J^{n-1}$ and $b\in J^n_{1}$ where $J_1=(x_2,\ldots,x_n).$ By \cite[Theorem 2.1]{Hu}, for all $n\geq 1,$ $$x_1(r-a)\in J^n_{1}\cap (x_1)\subset x_1J^{n-1}_{1}.$$ Therefore $x_1r\in x_1J^{n-1}.$ 
\eepf
\noindent
\begin{remark} Note that if $J=(x_1,\ldots,x_n)$ where $x_1,\ldots,x_n$ is a $d$-sequence and $\grade (J)=r$ then $x_1,\ldots,x_r$ is a regular sequence.\end{remark}
Next we state some well-known properties of  $d$-sequence. Here $G(J)=\oplus_{n\geq 0}J^n/J^{n+1}$ denotes the graded associated ring of $J.$
\bp\label{n3}
Let $(R,\mm)$ be a Noetherian local ring and $J$ an ideal in $R$ with $\grade (J)=s.$ Suppose $J$ is generated by a $d$-sequence $a_1,\ldots,a_s,\ldots,a_t.$ Then $\grade G(J)_{+}=\grade (J).$  
\ep

\bp\label{as}
Let $x_1,\ldots,x_s$ be a $d$-sequence in a Noetherian local ring $(R,\mm)$ and $J$ be an ideal minimally generated by $x_1,\ldots,x_s.$ Suppose $x_1$ is a nonzerodivisor on $R.$ Then $l(J/(x_1))=l(J)-1.$ 
\ep
The following theorem gives a lower bound for $\deg P(I/J)$ and using the lower bound of $\deg P(I/J)$ we show that if  $J$ is a complete intersection ideal then $\deg P(I/J)=l(J)-1.$
\bt\label{n4}
Let $(R,\mm)$ be a Noetherian local ring, $J\subsetneq I$ be ideals in $R$ such that  $\lm(I/J)<\infty.$ Then the following are true.
\ben
\item $\grade G(J)_{+}-1\leq \deg P(I/J).$
\item If $J$ is a reduction of $I$ then $\grade G(J)_{+}-1\leq \deg P(I/J)\leq l(J)-1.$ 
\item If $J$ is a complete intersection ideal and reduction of $I$ then $\deg P(I/J)=l(J)-1.$
\een
\et
\bpf
$(1)$ Let $S=R[X]_{\mm[x]}.$ Then $S$ is faithfully flat extension of $R,$ $S$ has infinite residue field, $\lm_S(I^nS/J^nS)=\lm_R(I^n/J^n)\lm_S(S/\mm S)=\lm_R(I^n/J^n)$ whenever $\lm_R(I^n/J^n)<\infty$ and $\grade G(JS)_{+}\geq \grade G(J)_{+}.$ Therefore without loss of generality we assume that $R$ has infinite residue field. Let $\grade G(J)_{+}=s.$ We may assume $s\geq 1.$ Then there exist $a^*_1,\ldots,a^*_s\in G(J)_{1}$ such that $a^*_1,\ldots,a^*_s$ is $G(J)$-sequence. Therefore for all $i=0,\ldots,s-1$ and $n\geq 0,$ we have $(J^{n+1}_i:(a^{(i)}_{i+1}))=J^{n}_i$ and hence $ (J^{n+1}_i:(a^{(i)}_{i+1}))\cap I^{n}_i=J^{n}_i$ where  $a^{(i)}_{i+1}$ denotes the image of $a_{i+1}$ in $R_i$ ($R_i, I_i, J_i$ are defined before Lemma \ref{need}). Thus by Lemma \ref{need}, for all $i=0,\ldots,s-2$ and $n\gg 0,$ we have \beqnn\label{eq1}\lm(I^{n+1}_{i+1}/J^{n+1}_{i+1})\leq \lm(I^{n+1}_i/J^{n+1}_i)-\lm(I^{n}_i/J^{n}_i).\eeqnn
Suppose $\deg P(I/J)<s-1.$ Then using the inequality (\ref{eq1}), for all $i=0,\ldots,s-2,$ we get that $\lm(I^{n}_{s-1}/J^{n}_{s-1})$ is a polynomial type function of degree less than zero and hence there exists an integer $k$ such that $I^{n}_{s-1}=J^{n}_{s-1}$ for all $n\geq k.$
\\We show that if $I^{m}_{s-1}=J^{m}_{s-1}$ for some integer $m$ then $I^{m-1}_{s-1}=J^{m-1}_{s-1}.$ Let $x'\in I^{m-1}_{s-1}$ where $x'$ denotes the image of $x$ in $R_{s-1}.$ Then $x'a^{(s-1)}_{s}\in I^{m}_{s-1}= J^m_{s-1}.$ Thus $x'\in J^m_{s-1}:a^{(s-1)}_{s}=J^{m-1}_{s-1}.$ Using this technique $m-1$ times, we get $I_{s-1}=J_{s-1}.$ This implies $I=J$ which is a contradiction. Therefore $\grade G(J)_{+}-1\leq \deg P(I/J).$
\\$(2)$ This follows from part $(1)$ and Remark \ref{n2}.
\\$(3)$ since  $J$ is a complete intersection ideal, $J$ is generated by a regular sequence and by Proposition \ref{n3}, we have $\grade G(J)_+=\grade(J)=\height J=l(J).$ Hence the result follows from part $(2).$
\eepf
\begin{example}
	{\em Let $R=K[X,Y]_{(X,Y)}$ where $K$ is a field. Let $J=(XY^2,X^4)$ and $I=(X^4,XY^2,X^3Y).$ Then $\lm(I/J)<\infty$ and $J$ is a reduction of $I.$ Now $\grade G(J)_{+}\geq 1,$ $J$ is generated by a $d$-sequence. For all $n\gg 0,$ $$I^n=(X^{4n-1}Y,X^{4n-4}Y^3,X^{4n-7}Y^5,\ldots, X^{4n-3n+2}Y^{2(n-1)+1})+J^n$$ where $J^n=(X^{4n},X^{4n-3}Y^2,X^{4n-6}Y^4,\ldots, X^{4n-3n}Y^{2n}).$ Then $\lm(I^n/J^n)$ is a polynomial in $n$ of degree one.}
\end{example}
The following theorem provides a better lower bound for $\deg P(I/J)$ and as a consequence of that we obtain a class of ideals for which the polynomial $P(I/J)(X)$ attains the maximal degree.
\bt\label{gra}
Let $(R,\mm)$ be a Noetherian local ring, $J\subsetneq I$ be ideals in $R$ such that $\lm(I/J)<\infty$ and $\grade (J)=s.$ Let $a_1,\ldots,a_s,b_1,\ldots,b_t$ be a $d$-sequence in $R$ which minimally generates the ideal $J$ and $t\geq 1.$  Suppose $I\cap(J':b_1)=J'$ where $J'=(a_1,\ldots,a_s)$ $(J'=(0)$ if $s=0)$. Then the following are true.
\ben
\item $\grade (J)\leq \deg P(I/J).$
\item If $J$ is a reduction of $I$ then $\grade (J)\leq \deg P(I/J)\leq \l(J)-1.$ 
\een
\et
\bpf 
$(1)$ ($R_i, I_i, J_i$ are defined before Lemma \ref{need}). 
\\First we show that if $s=0$ then $ H^0_{\R(J)_{+}}(\R(I))_1=0.$ Let $at\in H^0_{\R(J)_{+}}(\R(I))_1.$ Then for some $m>0,$ we have $({b^{m}_1}t^m)(at)=0$ in $\R(I).$ Hence $ab^{m}_1=0.$ If $m\geq 2,$ since $b_1,\ldots,b_t$ is a $d$-sequence in $R,$ we have $$ab^{m-2}_1\in ((0):b^2_1)=((0):b_1).$$ Thus $ab^{m-1}_{1}=0.$ Using the same method $m-1$-times, we get $ab_1=0.$ Therefore $a\in I\cap((0):b_1)=(0).$ 

\vspace*{2mm}
\noindent
Note that images of $a_i,\ldots,a_s,b_1,\ldots,b_t$ in $R_{i-1}$ is a $d$-sequence for all $i=1,\ldots,s.$ Now $\grade (J)=s$ implies $a_1,\ldots,a_s$ is a regular sequence in $R.$ Thus by Proposition \ref{n1}, for all $i=0,\ldots,s-1$ and $n\gg 0,$ we have $$ (J^{n+1}_i:(a^{(i)}_{i+1}))\cap I^{n}_i=J^{n}_i.$$ Hence by Lemma \ref{need}, for all $i=0,\ldots,s-1$ and $n\gg 0,$ we have \beqnn\label{eq2}\lm(I^{n+1}_{i+1}/J^{n+1}_{i+1})\leq \lm(I^{n+1}_i/J^{n+1}_i)-\lm(I^{n}_i/J^{n}_i).\eeqnn
Suppose $\deg P(I/J)\leq s-1.$ Then using the inequality (\ref{eq2}), for all $i=0,\ldots,s-1,$ we get that $\lm(I^{n}_{s}/J^{n}_{s})$ is a polynomial type function of degree less than zero and hence there exists an integer $k$ such that $I^{n}_{s}=J^{n}_{s}$ for all $n\geq k.$ 
\\ Now replace $R$ by $R_s$ and $I,J$ by $I_s,J_s$ respectively. Hence $J\subsetneq I,$ $J=(c_1,\ldots,c_t)$ where $c_1,\ldots,c_t$ is a $d$-sequence in $R$ $(c_1,\ldots,c_t\mbox{ are images of } b_1,\ldots,b_t\mbox{ in } R_s\mbox{ respectively}),$ $\grade(J)=0,$ $J'=(0)$ and $I\cap((0):c_1)=(0).$ Hence $ H^0_{\R(J)_{+}}(\R(I))_1=0.$ \\Consider the exact sequence of $\R(J)$-modules $$0\longrightarrow \R(J)\longrightarrow \R(I)\longrightarrow \R(I)/\R(J)\longrightarrow 0$$ which induces a long exact sequence of local cohomology modules whose $n$-graded component is $$\cdots\longrightarrow H^i_{\R(J)_{+}}(\R(I))_n\longrightarrow H^i_{\R(J)_{+}}(\R(I)/\R(J))_n\longrightarrow H^{i+1}_{\R(J)_{+}}(\R(J))_n\longrightarrow\cdots.$$ Consider the case $i=0$ and $n=1.$ Since $I^n=J^n$ for all $n\gg 0,$ $\R(I)/\R(J)$ is $\R(J)_{+}$-torsion. Thus $H^0_{\R(J)_{+}}(\R(I)/\R(J))_1=I/J.$ Since $J$ is generated by $d$-sequence, by Theorem \cite[Corollary 5.2]{Trung}, $\reg(\R(J))=0$ and hence $H^1_{\R(J)_{+}}(\R(J))_1=0.$ Therefore from the long exact sequence, for $i=0$ and $n=1,$ we get $I=J$ which is a contradiction. Hence $s\leq \deg P(I/J).$
\\$(2)$ Suppose $J$ is a reduction of $I.$ Then by part $(1)$ and Remark \ref{n2}, we get the required result.
\eepf
\begin{example}
	{\em Let $R=K[X,Y,Z,W]_{(X,Y,Z,W)}$ where $K$ is a field. Let $$I=(XZ,XW,YZ,YW)\mbox{  and }J=(XZ,YW,XW+YZ).$$ Then $J$ is a reduction of $I$ and $\lm(I/J)<\infty.$ By computations on Macaulay2 \cite{GS}, we get $\grade(J)=2,$ $I\cap((XZ,YW):(XW+YZ))=(XZ,YW)$ and $XZ,YW,XW+YZ$ is a $d$-sequence. Consider the ideals $\ov J =J/(XZ,YW)$ and $\ov I=I/(XZ,YW)$ in $\ov R=R/(XZ,YW).$ Let $x,y,z,w$ denote images of $X,Y,Z,W$ in $\ov R.$ Then $\ov{J}^n=(x^nw^n+y^nz^n)$ and $\ov{I}^n=(x^nw^n,y^nz^n).$ Hence $\lm(\ov{I}^n/\ov{J}^n)$ is nonzero constant for all $n\gg 0.$ Thus $2\leq \deg P(I/J).$ Since $l(J)=l(I)=\dim K[XZ,XW,YZ,YW]=3,$ we get $\deg P(I/J)=2.$}
\end{example}
Let $(R,\mm)$ be a Cohen-Macaulay local ring of dimension $d\geq 1.$ We say an ideal $I$ satisfies $G_s$ if $\mu(I_p) \leq \height p$ for every $p\in V(I)$ such that $\height(p)<s.$ A proper ideal $J$ is an {\it $s$-residual intersection} of $I$ if there exist $s$ elements $x_1,\ldots, x_s$ of $I$ such that $J = (x_1,\ldots,x_s) : I$ and $\height J \geq s.$ We say $J$ is a {\it geometric $s$-residual intersection} of $I$ if in addition we have $\height I + J> s.$ The ideal $I$ is said to satisfy the {\it Artin-Nagata property} $AN_{s}^-$ if the ring $R/K$ is Cohen-Macaulay for every $0\leq i\leq s$ and for every geometric $i$-residual intersection $K$ of $I.$
\begin{corollary}\label{Artin}
Let $(R,\mm)$ be a Cohen-Macaulay local ring of dimension $d\geq 2$ with infinite residue field, $I$ be an ideal of $R$  and $J$ be any  minimal reduction of $I$ with $\lm(I/J)<\infty.$  Suppose $I$ satisfies $G_{l(I)}$ and $AN_{l(I)-2}^{-}$ conditions. Then for any ideal $K$ with  $J\subsetneq K\subseteq I,$  $\deg P(K/J)=l(J)-1.$
\end{corollary}
\bpf
By \cite{Ul} and \cite{JU},  $\grade(J)\geq l(I)-1,$ $x_1,\ldots,x_{l(I)}$ is a $d$-sequence and $$K\cap((x_1,\ldots,x_{{l(I)}-1}):x_{l(I)})\subseteq I\cap((x_1,\ldots,x_{{l(I)}-1}):x_{l(I)})=(x_1,\ldots,x_{{l(I)}-1}).$$  If $\grade (J)=l(I)$ then by Theorem \ref{n4},  $\deg P(K/J)=l(J)-1.$ Suppose  $\grade(J)=l(I)-1.$ 
Then  by Theorem \ref{gra}, we obtain $\deg P(K/J)=l(J)-1.$
\eepf
\section{degree of the multiplicity function $e(I^n/J^n)$}
Let $(R,\mm)$ be a formally equidimensional local ring and $J\subsetneq I$ be ideals in $R.$ Then by \cite[Lemma 2.2]{cc15}, the chain of ideals $$\sqrt{J:I}\subseteq\cdots\subseteq\sqrt{J^n:I^n}\subseteq\cdots$$ stabilizes. Choose an integer $r\geq 1$ such that $K:=\sqrt{J^n:I^n}=\sqrt{J^{n+1}:I^{n+1}}$ for all $n\geq r.$ Put $\dim R/K=t.$ Then using associativity formula, for all $n\geq r,$ we have \beqnn\label{dd} e(I^n/J^n)=\sum\limits_{K\subseteq p, \dim R/p=t}e(R/p)\lm(I^n_p/J^n_p).\eeqnn 
By \cite[Proposition 2.4]{cc15}, $e(I^n/J^n)$ is a function of polynomial type of degree at most $\dim R-t.$ 
\\Throughout this section we use the following notation: $K=\sqrt{J^n:I^n}=\sqrt{J^{n+1}:I^{n+1}}$ for all $n\geq r,$ $t=\dim R/K$ and $S=\{p\in \spec R: K\subseteq p, \dim R/p=t \}.$ \\Next we show that if $J$ is a complete intersection ideal we can explicitly detect the $\deg e(I^n/J^n).$
\bp\label{d1}
Let $(R,\mm)$ be a formally equidimensional local ring, $J\subsetneq I$ be ideals in $R.$
\ben
\item Then $\grade G(J)_{+}-1\leq\deg e(I^n/J^n)\leq l(J).$
\item If $J$ is a reduction of $I$ then $\grade G(J)_+-1\leq\deg e(I^n/J^n)\leq l(J)-1.$
\item If $J$ is a complete intersection ideal then 
\ben
\item [$(a)$]$J$ is a reduction of $I$ implies $\deg e(I^n/J^n)=l(J)-1.$ 
\item [$(b)$]$J$ is not a reduction of $I$ implies $\deg e(I^n/J^n)=l(J).$ 
\een
\een
\ep
\bpf
$(1)$ ($K,S$ are defined at the beginning of the Section $3$). Suppose $J_p$ is reduction of $I_p$ for all $p\in S$ then by Theorem \ref{n4}, for all $p\in S,$  $$\grade G(J)_+-1\leq \grade G(J_p)_+-1\leq\deg P(I^n_p/J^n_p)\leq l(J_p)-1\leq l(J)-1.$$ Hence from the equation (\ref{dd}), $\grade G(J)_+-1\leq\deg e(I^n/J^n)\leq l(J)-1.$ 
\\Suppose $J_q$ is not a reduction of $I_q$ for some $q\in S.$ Therefore $J$ is not a reduction of $I$ and by \cite[Proposition 3.6.3]{FCV}, for all $n\gg0,$ we have $\height K\leq \height{\sqrt{JI^{n-1}:I^n}}\leq l(J).$ Then  by \cite[Theorem 2.1]{Rees}, $\height J\leq\deg P(I_q/J_q)=\height q\leq l(J).$ Therefore the result follows from the equation (\ref{dd}).
\\$(2)$ If $J$ is a reduction of $I$ then $J_p$ is reduction of $I_p$ for all $p\in S$ and the result follows from the argument in $(1).$
\\$(3)$ Since $J$ is a complete intersection ideal, we have $\grade(J)=l(J).$ Therefore by Proposition \ref{n3} and part $(2),$ we get that $J$ is a reduction of $I$ implies $\deg e(I^n/J^n)=l(J)-1.$ If $J$ is not a reduction of $I,$ then the result follows from \cite[Theorem 2.6]{cc16}. 
\eepf
\bl\label{impreq}
Let $(R,\mm)$ be a formally equidimensional local ring, $J\subsetneq I$ be ideals in $R.$ Suppose there exists a prime ideal $Q\supseteq J$ such that $\dim R/J=\dim R/Q$ and $I\nsubseteq Q.$ Then $\deg e(I^n/J^n)=\height J.$ In particular, if $J$ is unmixed (i.e. all associated primes of $J$ have same height) and $\sqrt{J}\subsetneq \sqrt{I},$ then $\deg e(I^n/J^n)=\height J.$
\el
\bpf
Let $\mbox{Min}(R/K)=\{Q_1,\ldots,Q_k\}$ ($K$ is defined at the beginning of the Section $3$). Put $T=Q_1\cdots Q_k.$ Since $K=\sqrt{J^r:I^r},$ there exists an integer $c>0$ such that $T^cI^r\subseteq J^r.$ Let $J=\bigcap\limits_{P\in Ass(R/J)}q(P)$ be an irredundant primary decomposition of $J.$ Then $$I^r\subseteq J:T^{\infty} \subseteq\bigcap\limits_{P\nsupseteq T,P\in Ass(R/J)}q(P).$$ Since $I\nsubseteq Q,$ we have $T\subseteq Q.$ Therefore $Q_j\subseteq Q$ for some $j$ and hence $Q=Q_j.$ Thus $J_{Q_j}$ is not a reduction of $I_{Q_j}=R_{Q_j}$ and $\dim R/Q=\dim R/K.$ Hence by \cite[Theorem 2.1]{Rees} and the equation (\ref{dd}), we get $\deg e(I^n/J^n)=\deg\lm(I^n_{Q_j}/J^n_{Q_j})=\dim R_{Q_j}=\height J.$
\eepf
We provide some sufficient conditions on $J$ such that for any $J\subsetneq I,$ we have $\sqrt{J:I}=\sqrt{J^n:I^n}$ for all $n\geq 1.$
\bp\label{req}
Let $(R,\mm)$ be a Noetherian local ring and $J\subsetneq I$ be ideals in $R$ such that one of the following conditions hold, 
\ben \item[(1)] the residue field of $R$ is infinite and $\grade G(J)_+\geq 1,$
 \item[(2)] $\grade (J)\geq 1$ and $J$ is generated by $d$-sequence, 
 \item[(3)] $J^{k+1}:J=J^k$ for all $k\geq 0.$\een Then $\sqrt{J:I}=\sqrt{J^n:I^n}$ for all $n\geq 1.$
\ep
\bpf
By \cite[Lemma 2.2]{cc15}, we know $\sqrt{J:I}\subseteq\sqrt{J^n:I^n}$ for all $n\geq 1.$ Suppose $x\in \sqrt{J^n:I^n}$ for any $n>1.$ Then $x^rI^n\subseteq J^n$ for some $r\geq 1.$ 
\\First consider the case that $\grade G(J)_+\geq 1$ and the residue field of $R$ is infinite. Then there exists an element $y\in J$ such that $J^n:(y)=J^{n-1}$ for all $n\geq 1.$ Then $x^ryI^{n-1}\subseteq x^rII^{n-1}\subseteq J^n$ implies $x^rI^{n-1}\subseteq J^n:(y)=J^{n-1}.$ Using this technique $n-1$ times we get $x^rI\subseteq J.$
\\ Note that in condition $(1),$ we require infinite residue field just to get a $G(J)$-regular element of degree $1.$ If $\grade (J)\geq 1$ and $J$ is generated by a $d$-sequence, say $J=(a_1,\ldots,a_s),$ by Propositions \ref{n1} and \ref{n3}, we get a $G(J)$-regular element of degree $1.$ Hence the result follows from the previous paragraph.
\\ If $(3)$ holds then for all $n\geq 1,$ $x^rJI^{n-1}\subseteq x^rII^{n-1}\subseteq J^n$ implies $x^rI^{n-1}\subseteq J^n:J=J^{n-1}.$
\eepf
The following theorem gives characterization of reduction in terms of $\deg e(I^n/J^n).$  
\bt\label{ad1}
Let $(R,\mm)$ be a formally equidimensional local ring  of dimension $d\geq 2,$ $J\subsetneq I$ be ideals in $R$ and $J$ has analytic deviation one. Suppose $l(J_p)< l(J)$ for all prime ideals $p$ in $R$ such that $\height p=l(J).$ Then the following are true.
\ben
\item If $J$ is not a reduction of $I$ then $\deg e(I^n/J^n)=l(J)-1.$
\item If $l(J)=d-1,$ $\depth (R/J)>0$ and for all $n\geq 1,$ $\sqrt{J:I}=\sqrt{J^n:I^n}$ then \\$J$ is a reduction of $I$ if and only if $\deg e(I^n/J^n)\leq l(J)-2.$
\een
\et
\bpf
(Note that the condition $l(J_p)< l(J)$ for all prime ideals $p$ in $R$ such that $\height p=l(J)$ implies $l(J)\leq d-1.$)
First we show that if $\height K=l(J)$ then $J$ is a reduction of $I$ (where $K=\sqrt{J^r:I^r},$ defined at the beginning of the section $3$).\wlg we may asumme that  $\ass(R/{\ov{J^n}})=\ass(R/{\ov{J^{n+1}}})$ for all $n\geq r.$ Let $q\in$Ass $(R/{\ov {J^r}}).$ We show that $(J^r:I^r)\nsubseteq q.$ If possible suppose $(J^r:I^r)\subseteq q.$ Then either $\height q>l(J)\geq l(J_q)=l(J^r_q)$ or $\height q= l(J)>l(J_q)=l(J^r_q).$ Therefore by \cite{MacA}, \cite[Lemma 3.1]{Rees}, $q\notin \ass (R/\ov{J^{nr}})=\ass (R/\ov{J^{r}})$ for any $n\geq 1,$ which is a contradiction. Thus by prime avoidance lemma, we have an element $$x\in (J^r:I^r)\setminus \underset{Q\in \ass (R/{\ov {J^r}})}\cup Q.$$ Therefore $x^{\prime}$ is a nonzerodivisor in $R/\ov{J^r}$ where $\prime$ denotes the image in $R/\ov{J^r},$ i.e. $\ov{J^r}:(x)=\ov{J^r}.$ This implies $I^r\subseteq J^r:(x)\subseteq \ov{J^r}:(x)=\ov{J^r}.$ Thus $J^r$ is a reduction of $I^r$ and hence $J$ is a reduction of $I.$
\\$(1)$ Since $J$ is not a reduction of $I,$ by \cite[Proposition 3.6.3]{FCV} and the previous paragraph, $\height K=\height J$ (where $K=\sqrt{J^r:I^r},$ defined at the beginning of the section $3$). Now $J$ is not a reduction of $I$ implies $J^r$ is not a reduction of $I^r.$ Therefore, by \cite{MacA},\cite{Rees5},\cite[Theorem 8.21]{Vas}, we have $J^r_p$ is not a reduction of $I^r_p$ for some prime ideal $p$ such that $K\subset p$ and $\height p=l(J_p).$
\\ Let $K\subset p\in \spec(R)$ such that $\height p\geq l(J)=\height K+1.$ Then either $\height p>l(J)\geq l(J_p)$ or $\height p= l(J)>l(J_p).$ Therefore $J^r_p$ is not a reduction of $I^r_p$ for some prime ideal $p\in S$ ($S$ is defined at the beginning of the section $3$). Hence $J_p$ is not a reduction of $I_p$ for some prime ideal $p\in S$ and thus by \cite[Theorem 2.1]{Rees}, $\deg e(I^n/J^n)=\height p=\height J=l(J)-1.$ 
\\$(2)$ It is enough to show that if $J$ is a reduction of $I$ then $\deg e(I^n/J^n)\leq l(J)-2.$
First we show that if $J$ is a reduction of $I$ then $\height K\leq l(J).$ If $\height K>l(J)$ then $K$ is $\mm$-primary and hence there exists an integer $n>0$ such that $\mm^nI\subseteq J.$ Since $\mm\notin\ass(R/J),$ we have  $I\subseteq J$ which is a contradiction. Suppose $\height K=\height J=l(J)-1.$ Since $J$ is a reduction of $I$, $I_p,J_p$ are $pR_p$-primary for all $p\in S$ and hence we get $\deg e(I^n/J^n)\leq l(J)-2.$ Now suppose $\height K=l(J).$ Since $J$ is a reduction of $I,$ for any $p\in S,$ we have $\height p=l(J)$ and hence $\deg P(I_p/J_p)\leq l(J_p)-1\leq l(J)-2.$ Thus $\deg e(I^n/J^n)\leq l(J)-2.$
\eepf
\begin{example}\label{ex2}
	{\em  We can not omit the condition $\depth(R/J)>0$ from Theorem \ref{ad1} $(2)$.  
		\\Let $R=\mathbb{Q}[x,y,z,w]_{(x,y,z,w)}$ and $J=(xz,yw,xw+yz).$ Then $3=l(J)=\height J+1,$ $\mm\in$ Ass$(R/J)$ and $l(J_p)\leq l(J)-1$ for all prime ideals $p$ in $R$ such that $\height p=l(J).$ Let $I=(xz,xw,yz,yw).$ Then $J$ is a reduction of $I$ and $\sqrt{J:I}=(x,y,z,w).$ Since $xz,yw,xw+yz$ is a $d$-sequence, by Proposition \ref{req}, $\sqrt{J:I}=\sqrt{J^n:I^n}$ for all $n\geq 1.$ Then $\deg e(I^n/J^n)=\deg \lm(I^n/J^n)=2.$}
\end{example}
\begin{example}\label{ex1}
	{\em This example shows existence of an ideal $J$ which satisfies conditions of Theorem \ref{ad1}.
		Let $R=\QQ[x,y,z,w]_{(x,y,z,w)}$ and $J=(xyw^2,xyz^2,xw^2+yz^2).$ Then by computations on Macaulay2 \cite{GS}, we get $l(J)=3,$ $J$ is an ideal of linear type (hence basic) and Ass$(R/J)=\{p_1=(w,z),p_2=(w,y),p_3=(x,z),p_4=(x,y)\}.$ Since $xyzw\in R\setminus J$ is integral over $J,$ $J$ is not integrally closed. Note that $J_{p_1}=(w^2,z^2)R_{p_1},J_{p_2}=(y,w^2)R_{p_2},J_{p_3}=(x,z^2)R_{p_3},J_{p_4}=(xy,xw^2+yz^2)R_{p_4}.$ Hence $J$ is generically complete intersection. Therefore by \cite[Lemma 2.5]{CZ}, $\grade G(J)_+\geq 1.$ Thus by Proposition \ref{req}, $\sqrt{J:I}=\sqrt{J^n:I^n}$ for all $n\geq 1$ and any ideal $I$ containing $J.$ Suppose $p$ is a prime ideal in $R$ of height $3.$ Then either $x\notin p \mbox{ (hence } J_p=(yz^2,w^2)R_p)$ or $y\notin p \mbox{ (hence } J_p=(xw^2,z^2)R_p)$ or $z\notin p \mbox{ (hence } J_p=(xy,xw^2+yz^2)R_p)$ or $w\notin p \mbox{ (hence } J_p=(xy,xw^2+yz^2)R_p).$ Thus $l(J_p)\leq l(J)-1.$ }
		\end{example}


\begin{thebibliography}{99}
	\bibitem{Am} J. Amao, {\em On a certain Hilbert polynomial}, J. London Math. Soc. (2) 14 (1976) 13-20.
	
	\bibitem{cc15} C. Ciuperc\u{a}, {\em Asymptotic growth of multiplicity functions}, J. Pure. Appl. Algebra, {\bf 219}, 2015, 1045-1054. 
	
	\bibitem{cc16} C. Ciuperc\u{a}, {\em Degrees of multiplicity functions for equimultiple ideals}, J. Algebra, {\bf 452}, 2016, 106-117.
	
	\bibitem{CZ}  T. Cortadellas, S. Zarzuela, {\em Burch's inequality and the depth of the blow up rings of an ideal}, J. Pure Appl. Algebra, {\bf 157}, 2001, 183-204.
	 
	 \bibitem{C2013} S. D. Cutkosky, Multiplicities associated to graded families of ideals,  Algebra and Number Theory 7 (2013), 2059 - 2083.
	 
	 \bibitem{C2014} S. D. Cutkosky, Asymptotic multiplicities of graded families of ideals and linear series,  Advances in Mathematics 264 (2014), 55 - 113.
	 
	 \bibitem{CSS} S. D. Cutkosky, P. Sarkar and H. Srinivasan, Mixed multiplicities of filtrations, to appear in Trans. Amer. Math. Soc., 2018. 
	 
	 	\bibitem{FCV} H. Flenner, L. O'Carroll and W. Vogel, {\em Joins and intersections}, Springer Monographs in Mathematics, Springer-Verlag, Berlin-Heidelberg-New York, 1999.
		
	\bibitem{GS} R. Grayson and M. E. Stillman, {EM Macaulay 2}, a software system for research in algebaric Geometry, Available at http://www.math.uiuc.edu/Macaulay2.
	
	 \bibitem{F} F. Hayasaka, Asymptotic vanishing of homogeneous components of multigraded modules and its applications, Journal of Algebra 513, 1-26, 2018. 
	
	\bibitem{HHRT} M. Herrmann, E. Hyry, J. Ribbe, Z. Tang, {\em Reduction numbers and multiplicities of multigraded structures}, J. Algebra,  {\bf 197}, 1997, 311-341.
	
	\bibitem{HPV} J. Herzog, T. J. Puthenpurakal and J. K. Verma {\em Hilbert polynomials and powers of ideals},  Mathematical Proceedings of the Cambridge Philosophical society, 2008, 623-642.
	
	\bibitem{Hu} C. Huneke, {\em The theory of d-sequences and powers of ideals}, Advances in Mathematics, 1982, 249-279.
		
	\bibitem{JU} M. Johnson and B. Ulrich, {\em Artin-Nagata properties and Cohen-Macaulay of associated graded rings}, Compos. Math., {\bf 103}, 1996, 7-29.
	
	\bibitem{KK} K. Kaveh and G. Khovanskii, Newton-Okounkov bodies, semigroups of integral points, graded algebras and intersection theory,  Annals of Math. 176 (2012), 925 - 978.
	
		\bibitem{KS} D. Kirby and D. Rees, Multiplicities in graded rings I : The general theory, Contemporary Math., 159 (1994), 209–267.
		
		\bibitem{LM} R. Lazarsfeld and M. Musta\c{t}\u{a}, Convex bodies associated to linear series, Ann. Sci. Ec. Norm. Super 42 (2009) 783 - 835.
	
	\bibitem{MacA} S. McAdam, {\em Asymptotic prime divisors}, Lecture Notes in Mathematics,
	{\bf 1023}, Springer, 1983.
	
	\bibitem{Ok} A. Okounkov, Why would multiplicities be log-concave?, in The orbit method in geometry and physics, Progr. Math. 213, 2003, 329-347.
	
	\bibitem{Rees4} D. Rees, {\em $\mathfrak a$-transforms of local rings and a theorem on multiplicities of ideals}, Proc. Cambridge Philos. Soc., {\bf 57}, 1961, 8-17.
	
	
	\bibitem{Rees2} D. Rees, {\em Generalizations of reductions and mixed multiplicities}, J. Lond. Math. Soc. (2) {\bf 29}, 1984, 397-414.
	
	
	\bibitem{Rees} D. Rees, {\em Amao's theorem and reduction criterion}, J. London Math. Soc., (2), {\bf 32}, 1985, 404-410.
	
	\bibitem{Rees5} D. Rees, {\em Reductions of modules}, Math. Proc. Camb. Phil. Soc., {\bf 101}, 1987, 431-449.
	
	\bibitem {P1} P. Samuel, {\em La notion de multiplicit\'e en alg\`ebre et en 
		g\'eom\'etrie alg\'ebrique}, J. Math. Pures Appl. {\bf 30 }(1951), 159-274.
	
	
	\bibitem{Trung} N. V. Trung, {\em The Castelnuovo regularity of the Rees algebra and the associated graded ring}, Trans. Amer. Math. Soc., 350, 1998, 2813-2832.
	
	\bibitem{Ul} B. Ulrich, {\em Artin-Nagata properties and reduction of ideals}, Contemp. Math., {\bf 159}, 1994, 99-113.
		
	\bibitem{Vas}  W. Vasconcelos, {\em Integral closure. Rees algebras, multiplicities, algorithms}, Springer Monographs in Mathematics, Springer-Verlag, Berlin, 2005.
	
	
\end{thebibliography}
	\end{document}